\documentclass{article}

\usepackage{amsmath}
\usepackage{amssymb}
\usepackage{amsfonts}
\usepackage{amsthm}
\usepackage{tikz}
\usepackage{color}
\usepackage{graphicx}
\usepackage{mathrsfs}
\usepackage{url}
\usepackage{setspace}
\usepackage{comment}

\newtheorem{thm}{Theorem}

\newtheorem{cnj}[thm]{Conjecture}

\newtheorem{fct}[thm]{Fact}

\def\a{{\alpha}}

\def\c{{\gamma}}

\def\d{{\delta}}
\def\D{{\Delta}}

\def\m{{\mu}}

\def\p{{\pi}}

\def\s{{\sigma}}
\def\t{{\tau}}

\def\x{{\chi}}

\def\bz{{\mathbf{0}}}

\def\bmf{{\mathbf{f}}}
\def\bmF{{\mathbf{F}}}

\def\cA{{\cal A}}
\def\cB{{\cal B}}
\def\cC{{\cal C}}
\def\cD{{\cal D}}

\def\cF{{\cal F}}

\def\cH{{\cal H}}
\def\cI{{\cal I}}

\def\cS{{\cal S}}

\def\xh{{\hat{\x}}}

\def\od{{\overline{d}}}
\def\oF{{\overline{F}}}
\def\oG{{\overline{G}}}

\def\oS{{\overline{S}}}

\def\deg{{\sf deg}}
\def\dist{{\sf dist}}
\def\ekr{{\sf EKR}}

\def\hk{{\sf HK}}

\def\kth{{k^{\rm th}}}
\def\knr{{K(n,r)}}

\def\da{{\downarrow}}

\def\mt{{\emptyset}}

\def\rar{{\rightarrow}}
\def\sse{{\subseteq}}

\def\sqr#1#2{{\vcenter{\hrule height.#2pt
        \hbox{\vrule width.#2pt height#1pt \kern#1pt
                \vrule width.#2pt}
        \hrule height.#2pt}}}
\def\gbox{{\mathchoice\sqr45\sqr45\sqr{2.1}3\sqr{1.5}3}}

\definecolor{brwn}{RGB}{140, 70, 20}
\definecolor{gren}{RGB}{  0,140, 10}

\author{
Glenn Hurlbert\thanks{
Department of Mathematics and Applied Mathematics,
Virginia Commonwealth University, 
Richmond, VA, USA, 
\texttt{ghurlbert@vcu.edu}.}
}

\title{A Survey of the Holroyd-Talbot Conjecture}

\date{}

\begin{document}

\maketitle

\begin{abstract}
A family of sets is \textit{intersecting} if every pair of its members has an element in common.
Such a family of sets is called a \textit{star} if some element is in every set of the family.
Given a graph $G$, let $\m(G)$ denote the size of the smallest maximal independent set of $G$.
In 2005, Holroyd and Talbot conjectured the following generalization of the Erd\H os-Ko-Rado Theorem: for $1\le r\le \m(G)/2$, there is a maximum size intersecting family of independent $r$-sets that is a star.
In this paper we present the history of this conjecture and survey the results that have supported it over the last 20 years.
\end{abstract}


\section{Introduction}
\label{s:Intro}

For $0\le r\le n$, set $[n]=\{1,2,\ldots,n\}$, $2^{[n]}=\{S\ \sse\ [n]\}$, and $\binom{[n]}{r} = \{S\ \sse\ [n]\mid |S|=r\}$. 
For $\cF\ \sse\ 2^{[n]}$, define $\cap\cF=\cap_{S\in\cF}S$, $\cF^r=\{S\in\cF\ \mid\ |S|=r\} = \cF\cap\binom{[n]}{r}$, and $\cF_X=\{S\in\cF\mid X\ \sse\ S\}$.
We use the notational conventions $A+x=A\cup\{x\}$ and $A-x=A\setminus\{x\}$.
For $1\le t\le r$, we call $\cF$ $t$-{\it intersecting} if $|A\cap B|\ge t$ for every $A,B\in\cF$.
For convenience, we write $\cF_x$ instead of $\cF_{\{x\}}$, and say {\it intersecting} in place of $1$-intersecting.
If $\cap\cF\not=\emptyset$, we say that $\cF$ is a {\it star} (an $r$-{\it star} if $\cF\ \sse\ \binom{[n]}{r}$); in this case, any $x\in\cap\cF$ is called a ${\it center}$.
If $x$ is the center of a maximum size star of $\cF$, we call $x$ a {\it maximum} center.
If $|\cap\cF|=k$ we say that $\cF$ is $k$-{\it centered}.
Erd\H{o}s, Ko, and Rado \cite{ErdKoRad} proved the following foundational result in extremal set theory, that, for large enough $n$, the biggest $t$-intersecting families of $r$-subsets of $[n]$ are $t$-centered $r$-stars.

\begin{thm}
\label{t:EKR}
{\bf (Erd\H os-Ko-Rado, 1961)}
For all $1\le t\le r$, there is an integer $M(r,t)$ such that if $n>M(r,t)$ and $\cF\subseteq \binom{[n]}{r}$ is $t$-intersecting then $|\cF|\leq \binom{n-t}{r-t}$, with equality if and only if $\cF=\binom{[n]}{r}_X$ for some $X\in\binom{[n]}{t}$.
\end{thm}

In 1978, Frankl \cite{FranklEKRtrue} proved the optimal value of $M(r,t)=(r-t+1)(t+1)$ for all $t\ge 15$, and Wilson \cite{Wilson} in 1984 used Delsarte's linear programming bound \cite{Delsarte} to prove the same for the remaining values of $t\le 14$.
When $n=(r-t+1)(t+1)$, the Erd\H{o}s-Ko-Rado bound still holds, but there are non-stars of maximum size as well.
For example, when $t=1$ we have $M(r,1)=2r$ and any family that includes exactly one of each complementary pair $\{S,\oS\}$ is intersecting and of maximum size.
When $t>1$, the Complete Intersection Theorem of Ahlswede-Khachatrian \cite{AhlsKhac} identifies the maximum size and structure of $t$-intersecting families in the range of $2r-t<n<2r$.
This result had been conjectured in \cite{FranklEKRtrue}, and significant progress had been made on it in \cite{FranFureBeyond}.
In this paper we will restrict our attention to the case $t=1$ only.

Since the publication of Theorem \ref{t:EKR}, in which the authors used shifting to achieve the result, a number of different methods have been used to prove it: 
bounding the size of the largest intersecting non-star family via shifting \cite{HiltMiln}; 
finding local bounds that extend to global bounds via the cycle method \cite{Katona}; 
applying the Kruskal-Katona Theorem to shadows of complements \cite{Daykin,KatonaShadInt}; 
providing an upper bound on the Shannon capacity of the Kneser graph \cite{Lovasz}; 
using Delsarte's linear programming bound \cite{Wilson}; applying the Hoffman-Delsarte ratio bound to the independence number of the Kneser graph \cite{GodsMeag}; 
using the method of linearly independent polynomials \cite{FurHwaWei}; 
displaying an injection from any intersecting family into a star \cite{FranFure,HurlKamaInject};
and using the Borel Fixed Point Theorem from algebraic group theory \cite{WoodroofeAlg}.

For a graph $G$, we write $\cI(G)$ for the family of all independent sets of $G$, and use $\cI^r(G)$ in place of $\cI(G)^r$ and $\cI_v(G)$ in place of $\cI(G)_v$, etc.
Suppose that $\cF\subseteq\cI^r(G)$ is an intersecting subfamily of maximum size.
Define $G$ to be $r$-\ekr\ if some $v\in V(G)$ satisfies $|\cI^r_v(G)|=|\cF|$, and {\it strictly} $r$-\ekr\ if every such $\cF$ equals $\cI_v^r(G)$ for some $v$.
Such a $v$ is called an $r$-{\it witness} of $G$ (either strictly so or not).
Let $\a(G)$ denote the independence number of $G$, and $\m(G)$ denote the size of a smallest maximal independent set of $G$; i.e., the size of the smallest independent dominating set of $G$.
Holroyd and Talbot \cite{HolrTalb} made the following conjecture to generalize the $t=1$ case of Theorem \ref{t:EKR}, which is the empty graph case of this conjecture.

\begin{cnj}
\label{j:HolTal}
{\bf (Holroyd-Talbot, 2005)}
For any graph $G$, if $1\leq r\leq \m(G)/2$, then $G$ is $r$-\ekr, and strictly so if $1<r<\m(G)/2$.
\end{cnj}

The necessity of the $\m(G)/2$ condition in this conjecture stems from their example of two disjoint%
\footnote{When we speak of disjoint graphs, we mean vertex-disjoint.} 
copies of the complete bipartite graph $K_{3,3}$: $G=2K_{3,3}$, with $E(G)=(A\times B)\cup (C\times D)$.
Here we have $\m(G)=6$ and a maximum star in $\cI^4(G)$ has size 20, while the family consisting of every independent 4-set with at least two elements of $A$ is an intersecting non-star of size 24.
Hence $G$ is not $4$-\ekr.
More generally, two disjoint copies of $K_{2k-1,2k-1}$ has $\m/2=2k-1$ and is not $2k$-\ekr.

Because every graph $G$ is technically $r$-\ekr\ (vacuously) when $r\le 0$ or $r>\a(G)$, we will drop the obvious condition ``$1\le r\le \a(G)$'' that would normally follow ``for all $r$'', only writing constraints when they are more restrictive.


\section{Precursors to the Holroyd-Talbot Conjecture}
\label{s:Precursors}

In 1972, Berge \cite{Berge1} considered intersecting subfamilies of the collection of $k$-cliques in $K_{n_1,\ldots,n_k}$, and proved the following.

\begin{thm}
\label{t:Berge}
For $k\ge 2$, if $\cF$ is an intersecting family of $k$-cliques in $G=K_{n_1,\ldots,n_k}$, with $n_1\le n_2\le\cdots\le n_k$, then $|\cF|\le n_2\cdots n_k$.
\end{thm}

Of course, if $\cF$ is a maximum intersecting family and $f(G)$ is the maximum number of $k$-cliques containing the same vertex of $G$, then $n_2\cdots n_k = f(G) \le |\cF| \le \x'(G)$, where $\x'$ denotes the {\it chromatic index}.
Berge achieved the result by proving that $\x'(G)=f(G)$.
The paper was also published in \cite{Berge2}.
The result can be considered the first precursor to the Holroyd-Talbot Conjecture because it shows that the complement $\oG$ is $k$-\ekr.
For the case that each $n_i=n>2$, Livingston \cite{Livingston} proved that $\oG$ is strictly $k$-\ekr.
Other proofs of this result were given later in \cite{Gronau,Moon}.

Deza and Frankl \cite{DezaFran} considered the following scenario.
An $\ell$-{\it signed set} is a pair $(A,f)$, where $A\ \sse\ [n]$ and $f:A\rar [\ell]$, and its {\it size} equals $|A|$.
For two such pairs $(A,f)$ and $(B,g)$ define $(A,f)\cap (B,g) = (C,h)$, where $C=\{x\in A\cap B\mid f(x)=g(x)\}$ and $h:C\rar [\ell]$ such that $h(x)=f(x)$. 
Say that a family of $\ell$-signed sets is $t$-{\it intersecting} if the size of the intersection of any two of its members is at least $t$.
They proved the following, first stated by Meyer in \cite{Meyer} (stated in the complementary form of Berge/Livingston).

\begin{thm}
\label{t:Signed}
Let $r\le n$, $\ell\ge 2$, and $\cF\ \sse\ \binom{[n]}{r}$ be a $t$-intersecting family of $\ell$-signed sets.
Then $|\cF|\le \binom{n-1}{r-1}\ell^{r-1}$.
\end{thm}

{
\color{brwn}
\noindent
\hrulefill
\medskip

\centerline{
***
$\quad$
Theorem \ref{t:Signed} is stated incorrectly.
It should really be the following theorem.
Set $n_0(r,1)=r$.
$\quad$
***
}
\medskip

{\it 
\noindent
{\bf Theorem 4.}
Let $n\ge n_0(r,t)$, $\ell\ge 2$, and $\cF\ \sse\ \binom{[n]}{r}$ be a $t$-intersecting family of $\ell$-signed sets.
Then $|\cF|\le \binom{n-t}{r-t}\ell^{r-t}$.
}
\medskip

\noindent
\hrulefill
}
\medskip

When $t=1$, $\ell=2$, and $r\le n/2$, the result follows from Theorem \ref{t:EKR} and the following Fact.

\begin{fct}
\label{f:Half}
If $\cF\ \sse\ 2^{[n]}$ is intersecting, then $|\cF|\le 2^{n-1}$.
\end{fct}

Fact \ref{f:Half} holds because at most one of each pair $\{F,\oF\}$ can be a member of $\cF$, and is best possible because of $2^{[n]}_x$.
Thus, if $\cF$ is an intersecting family of 2-signed sets, then there are at most $\binom{n-1}{r-1}$ $r$-sets $A$ in some pair of $\cF$ when $r\le n/2$ by Theorem \ref{t:EKR}.
In addition, for each such $A$, Fact \ref{f:Half} implies that there are at most $2^{r-1}$ functions $f$ such that $(A,f)\in\cF$, since the set of those $f$ is isomorphic to the set of characteristic vectors of an intersecting family of subsets of $A$.
That is, for fixed $A\in\binom{[n]}{r}$, the set of pairs $(A,f)$, where $f:A\rar [2]$ is isomorphic to $2^A$ by the correspondence $B(A,f)=\{a\in A\mid f(a)=2\}$; clearly $(A,f)\cap (A',f')\not=\mt$ if and only if $B(A,f)\cap B(A',f')\not=\mt$.
But Theorem \ref{t:Signed} needs to be verified for the wider range $r\le n$, as well as for $t>1$ and for $\ell>2$.
Deza and Frankl \cite{DezaFran} prove this by using the shifting technique introduced in \cite{ErdKoRad}, which we define now.

For $1\le i<j\le n$, and for any family $\cF\ \sse\ 2^{[n]}$, define the shifting function $\s_{j,i}:\cF\rar 2^{[n]}$ by 
\[
\s_{j,i}(A) = \left\{
\begin{array}{ll}
A+i-j & \mbox{if}\ j\in A,\ i\not\in A,\ \mbox{and}\ A+i-j\not\in\cF,\ \mbox{and}\\
\s_{j,i}(A)=A & \mbox{otherwise.}
\end{array}
\right.
\]

We extend the shifting function to families of sets by $\s_{j,i}(\cF) = \{\s_{j,i}(A)\mid A\in\cF\}$.
It is well known that $|\s_{j,i}(\cF)| = |\cF|$ and that $\s_{j,i}(\cF)$ is $t$-intersecting whenever $\cF$ is.
Thus, when proving EKR-type results, one may always assume that an intersecting family $\cF$ is {\it shifted}; that is, $\s_{j,i}(\cF)=\cF$ for all $1\le i<j\le n$.
(See \cite{FranklShifting} for the many other wonderful properties and implications of shifting.)
The shifting technique as been applied to many varieties of intersection theorems, as well as to other types of problems in extremal set theory and other areas of combinatorics and graph theory.
Its power and utility makes it a primary tool in tackling Conjecture \ref{j:HolTal}.

The reason that Theorem \ref{t:Signed} is a precursor to Conjecture \ref{j:HolTal} is that the family of all $\ell$-signed $r$-subsets of $[n]$ corresponds to $\cI^r(G)$, where $G$ is the disjoint union of $n$ copies of the complete graph $K_\ell$, which we denote $nK_\ell$.
With this correspondence, we may think of the $t=1$ case of Theorem \ref{t:Signed} as being the second verification of Conjecture \ref{j:HolTal}.
In fact, this result goes farther than Holroyd-Talbot because it holds true for $r\le n=\m(nK_\ell)$.
Deza and Frankl \cite{DezaFran} also showed that $nK_\ell$ is strictly $r$-\ekr\ unless $\ell=2$ and $r=n\ge 3$.

Bollob\'as and Leader \cite{BollLead} gave an alternate proof of Theorem \ref{t:Signed} by using the idea underlying Katona's {\it cycle method} that proved Theorem \ref{t:EKR}.
The essence of the method is to define, for any permutation $\p$ of $[n]$, written cyclically, the family $\cF_\p$ of all sets of $\cF$ that appear consecutively on the cycle.
For any $A\in\cF_\p$, one can cut the cycle between any of the $r-1$ consecutive pairs of $A$, creating a pair of disjoint, consecutive sets on the cycle that begin at the adjacent elements on opposite sides of the cut and progress in opposite directions.
Since $\cF$ is intersecting, at most one of each pair is in $\cF_\p$, implying that $|\cF_\pi|\le r$.
Because each $A\in\cF$ is in a uniform number of $\cF_\p$, double-counting the pairs $(A,\cF_\p)$ with $A\in\cF_\p$ leads to the result.

In the language of independent sets in $nK_\ell$ instead of signed sets, a the $j^{\rm th}$ vertex $(i,j)$ in the $i^{\rm th}$ (maximum) clique of  $nK_\ell$ represents the element $x$ with sign $f(x)=i$, 
What Bollob\'as and Leader were able to do is make the double-counting work by only considering what they called {\it good} permutations of $[n\ell]$: permutations $\p$ in which, for each clique, the vertices of that clique are equally spaced around the cycle. 
This guarantees that consecutive $r$-sets in $\p$ are independent in $nK_\ell$, and that the necessary uniformity condition holds.

Recently, Feghali \cite{Feghali} gave an injective proof of the $t=1$ and $r\le n/2$ instance of Theorem \ref{t:Signed}, along the lines of \cite{FranFure,HurlKamaInject}.

The above results are for uniform families; in \cite{Bey}, Bey generalized the $t=1$ instance to antichains in general, proving the following theorem.
For an antichain $\cF\ \sse\ 2^{[n]}$, define is $f$-{\it vector} $\bmf(\cF) = (\bmf(r))_{r=0}^n$ by $\bmf(r) = |\cF^r|$.
For $0\le i\le n$, let $\bmf_i$ denote the $f$-vector with $\bmf_i(i)=\binom{n-1}{i-1}\ell^{i-1}$ and $\bmf_i(r)=0$ otherwise; that is, the $f$-vector of a star in $\cI^i(\oG)$, where $G=K_{n_1,\ldots,n_k}$.
Finally, set $\bz=(0)_{r=0}^n$ and write $\bmF$ for the convex hull of $\{\bz,\bmf_1,\ldots,\bmf_n\}$.

\begin{thm}
\label{t:Bey}
If $\cF\ \sse\ \cI(\oG)$ is an intersecting antichain, for $G=K_{n_1,\ldots,n_k}$ with $k\ge 2$, then $\bmf(\cF)\in\bmF$.
\end{thm}

For the instance of $t=1$ and $\ell=2$, Fuentes and Kamat \cite{FuenKama} proved the following common generalization of Theorems \ref{t:EKR} and \ref{t:Signed}.
Let $\cH^{(p,s)}(n)$ denote the family of vertex sets of subgraphs of $nK_2$ that have $p$ edges and $s$ isolated vertices.

\begin{thm}
\label{t:FuenKama}
Let $n\ge 2p+s$.
Then
\begin{enumerate}
    \item 
    if $n\ge 2p+2s$, then $\cH^{(p,s)}(n)$ is $(2p+s)$-\ekr, and strictly so when $n>2p+2s$,
    \item 
    $\cH^{(p,s)}(2p+s)$ is $(2p+s)$-\ekr, but not strictly so, and
    \item 
    if $n>2p+1$, then $\cH^{(p,1)}(n)$ is strictly $(2p+1)$-\ekr.
\end{enumerate}
\end{thm}

They conjecture that $\cH^{(p,s)}(n)$ is $(2p+s)$-\ekr\ for all $n\ge 2p+s$.

Another scenario involves families of subsets of integers, in which no subset contains consecutive elements.
Such subsets of $[n]$ correspond to independent sets in the path $P_n$.
Disallowing the additional pair $\{1,n\}$ corresponds to independent sets in the cycle $C_n$; such subsets are known as {\it separated} sets.
More generally, for a graph $G$, denote by $G^{(k)}$ the $\kth$ {\it power of} $G$: the graph on $V(G)$ having an edge $xy$ precisely when the distance in $G$ between $x$ and $y$ is positive and at most $k$.
Then a $k$-{\it separated} set is any set in $\cI(C_n^{(k)})$.
Talbot \cite{Talbot} proved the following result, stated in terms of $k$-separated sets, which had been conjectured by Holroyd and Johnson in \cite{HolrJohn} (originally posed by Holroyd \cite{Holroyd} in terms of sorties of King Arthur's knights sitting far apart from one another at the Round Table).

\begin{thm}
\label{t:Cycles}
The cycle power $C_n^{(k)}$ is $r$-\ekr\ for all $r$, and strictly so unless $k=1$ and $n=2r+2$.
\end{thm}

Talbot's proof uses induction, a simplified shift operator, a very careful partition of $\cF$ into $k+3$ parts, and seven pages of detailed analysis.
When $k=1$ and $n=2r+2$, a 1-centered $r$-star in $C_n^{(k)}$ has size $|\cI^{r-1}(P_{2r-1})|=\binom{r+1}{r-1}=\binom{r+1}{2}$, and an example of a non-star intersecting family of that size is given by the following construction.
Call a vertex of $V=V(C_{2r+2})$ {\it even} ({\it odd}) if its index is even (resp. odd).
For any $1\le i\le r$, any set $S$ of $i$ consecutive even vertices, there is a unique set $T=T(S)$ of $r-i$ odd vertices such that $F_S=S\cup T\in\cI^r(C_{2r+2})$.
Let $\cS$ denote the family of such $S$ for $i>(r+1)/2$; if $r$ is odd, also include half of the sets $S$ for $i=(r+1)/2$ by choosing one from each pair $\{S,\oS\}$, where the complement is over even vertices.
Then $\cS$ is intersecting of size $\binom{r+1}{2}$; hence $\cF=\{F_S\mid S\in\cS\}$ is intersecting of that size.

The {\it Kneser graph} $\knr$ consists of vertex set $\binom{[n]}{r}$, with edges between disjoint pairs of vertices.
Thus, intersecting families of $r$-subsets of $[n]$ correspond to independent sets in $\knr$, and Theorem \ref{t:EKR} is equivalent to the fact that $\a(\knr)=\binom{n-1}{r-1}$, established by Godsil and Meagher \cite{GodsMeag}.
For an interesting connection, Lov\'asz \cite{Lovasz} proved that the chromatic number of $\knr$ is $n-2r+2$ and, subsequently, Schrijver \cite{Schrijver} showed that the subgraph of $\knr$ induced by the separated sets is a minimal subgraph having this chromatic number. 

Building on Talbot's ideas, Holroyd, Spencer, and Talbot \cite{HolSpeTal} proved the following three analogous results for disjoint unions of complete graphs, for powers of paths, and for certain other disjoint unions of graphs.
This was the first paper to introduce the ``$r$-\ekr'' terminology.

\begin{thm}
\label{t:Cliques}
If $G$ is a disjoint union of $n\ge r$ complete graphs, each of size at least two, then $G$ is $r$-\ekr.
\end{thm}

Note that Theorem \ref{t:Cliques} is a common generalization of Theorem \ref{t:Berge} (to $r<k$) and the $t=1$ case of Theorem \ref{t:Signed} (to non-uniform part sizes).

\begin{thm}
\label{t:Paths}
The path power $P_n^{(k)}$ is $r$-\ekr\ for all $r$, and strictly so unless $k=1$ and $n=2r+1$.
\end{thm}

Theorem \ref{t:Paths} follows from Theorem \ref{t:Cycles} by a straightforward lemma that they prove: If $G$ is $r$-\ekr\ with center $v$ then, for any set $S\ \sse\ N_G(v)$ of neighbors of of $v$, $G-S$ is also $r$-\ekr\ with center $v$.
The trick in applying the lemma is to remove, not all neighbors of $v$ but, rather, just those on the clockwise side of $v$.

\begin{thm}
\label{t:Unions}
If $G$ is a disjoint union of $n\ge 2r$ complete graphs, cycles, and paths, including an isolated vertex ($K_1$ component), then $G$ is $r$-\ekr.
\end{thm}

Note that all three theorems verify Conjecture \ref{j:HolTal} for their given instances, while only Theorem \ref{t:Unions} is restricted to the range $r\le\m(G)/2$; it was later proven for the complete range of $r$ by Borg and Holroyd in \cite{BorgHolrSingle}.

Talbot's proof of Theorem \ref{t:Cycles} used a shift operator that is equivalent to contracting an edge of a graph.
This connection was made explicit in \cite{HolSpeTal}, where Holroyd, Spencer, and Talbot denoted by $G/e$ the graph obtained from $G$ by contracting the edge $e$.
Additionally, they defined $G\da e$ to be the graph obtained from $G$ by deleting the vertices of $e$ and their neighbors; i.e., $G-N(u)-N(v)$, where $e=uv$ and $N(x)$ denotes the neighbors of a vertex $x$.
They were able to show that if $\cF\ \sse\ \cI^r(G)$ is intersecting, $e=uv$ is an edge of $G$, and $x\in V(G\da e)$, then there are families $\cA$, $\cB$, $\cC$, and $\cD$ such that
\begin{itemize}
    \item 
    $|\cF|=|\cA|+|\cB|+|\cC|+|\cD|$,
    \item 
    $\cA\ \sse\ \cI^r(G/e)$ is intersecting,
    \item 
    $\cB\ \sse\ \cI^{r-1}(G\da e)$ is intersecting,
    \item 
    $\cC=\{F\in\cF_u\mid N(v)\cap (F-\{u\})\not=\mt\}$,
    \item 
    $\cD=\{F\in\cF_v\mid N(u)\cap (F-\{v\})\not=\mt\}$,
    \item 
    $B\in\cB$ and $E\in\cC\cup\cD$ implies that $B\cap E\cap V(G\da e)\not=\mt$,
    \item 
    $C\in\cC$ and $D\in\cD$ implies that $C\cap D\cap V(G\da e)\not=\mt$, and
    \item 
    $|\cI^r_x(G)| = |\cI^r_x(G/e)| + |\cI^{r-1}_x(G\da e)| + |\cC_x| + |\cD_x|$.
\end{itemize}
These properties are then used to prove the above results by induction on both $r$ and $n$.
For example, in the case of Theorem \ref{t:Cliques}, one chooses $u$ and $v$ from the same component to obtain that $\cC=\cD=\mt$, so that
\[
|\cF| = |\cA|+|\cB| \le |\cI^r_x(G/e)|+|\cI^{r-1}_x(G\da e)| = |\cI^r_x(G)|.
\]
Variations on this theme produce the other results, in some cases by choosing an appropriate vertex $x$.
In the case of Theorem \ref{t:Unions}, $x$ is an {\it isolated} vertex (a vertex with no neighbors), which makes sense by virtue of it being a member every maximal independent set.
Thus, importantly, knowing which vertices are maximum centers is sometimes necessary to show that a particular graph is \ekr.


\section{Further Support for the Holroyd-Talbot Conjecture}
\label{s:Support}

Conjecture \ref{j:HolTal} was made in \cite{HolrTalb}.
In that paper, Holroyd and Talbot proved the following theorem for the $r=2$ case.

\begin{thm}
\label{t:Alpha}
Let $G$ be an $n$-vertex graph of minimum degree $\d$.
\begin{enumerate}
    \item 
    If $\a(G)=2$, then $G$ is strictly $2$-\ekr.
    \item 
    If $\a(G)\ge 3$, then $G$ is $2$-\ekr\ if and only if $\d\le n-4$, and strictly so if and only if $\d\le n-5$, with any vertex of minimum degree as a maximum center.
\end{enumerate}
\end{thm}

They also noted an interesting set of natural examples: If $G$ is the graph of vertices and edges of a platonic solid $S$ then $G$ is $\a(G)$-\ekr\ if and only if $S$ is not the dodecahedron.

Additionally, they proved an \ekr\ result involving graph products.
For graphs $G$ and $H$, define the {\it lexicographic product} $G[H]$ to have vertex set $V(G)\times V(H)$, with edges $(u,x)(v,y)$ if $uv$ is an edge of $G$ or if $u=v$ and $xy$ is an edge of $H$.
For example, if $G=K_2$ and $H=mK_1$ (the empty graph on $m$ vertices) then $G[H]=K_{m,m}$ (the $m$-regular complete bipartite graph), while $H[G]=mK_2$ (the perfect matching on $2m$ vertices).

\begin{thm}
\label{t:LexProd}
If $G$ is $r$-\ekr\ and $m\ge 1$ then $G[K_m]$ is $r$-\ekr.
\end{thm}

They achieved this result by (a) exhibiting conditions that produce an $r$-witness and (b) using those conditions to show that if $v$ is an $r$-witness of $G$ then $(v,x)$ is an $r$-witness of $G[K_m]$ for all $x\in V(K_m)$. 

Furthermore, they generalized the lexicographic product as follows.
If $G$ has $n$ vertices $v_1,\ldots,v_n$ and $\cH=(H_1,\ldots,H_n)$ is a sequences of $n$ graphs then the graph $G[\cH]$ has vertices $\cup_{i=1}^n (\{v_i\}\times V(H_i))$, with edges $(v_i,x)(v_j,y)$ if $v_iv_j$ is an edge of $G$ or if $i=j$ and $xy$ is an edge of $H_i$.
Then $G[H]$ is the case in which each $H_i=H$.
They used this general construction to provide several informative examples, such as (1) that if $\cH$ is a sequence of complete graphs then the conclusion of Theorem \ref{t:LexProd} does not necessarily hold if the complete graphs are not all of the same size, and (2) that, for any $r\ge 3$, there are graphs of arbitrarily large independence number that are not $r$-\ekr.

Finally, they used shifting to prove the following theorem.

\begin{thm}
\label{t:TwoCompBip}
Let $G$ be a disjoint union of two complete multipartite graphs.
Then $G$ is $r$-\ekr\ for all $r\le \m(G)/2$ and strictly so unless $r=\m(G)/2$.
\end{thm}

Borg and Holroyd \cite{BorgHolrDouble} then proved much more.

\begin{thm}
\label{t:ManyCompBip}
Let $G$ be a disjoint union of complete multipartite graphs, at least one of which is a single vertex.
Then $G$ is $r$-\ekr\ for all $r\le \m(G)/2$ and strictly so unless $r=\m(G)/2$.
\end{thm}

Their work was phrased in the form of {\it double partitions}: a collection of disjoint {\it large} sets $A_i$, each of which is partitioned into {\it small} sets $A_{i,j}$, with every intersecting family of sets chosen from those that intersect each large set in some portion of exactly one of its small sets, if at all.

In 2009, Borg \cite{BorgExt} proved one of the more significant results to date, roughly that every graph $G$ is $r$-\ekr\ for every $r$ up to about $(\frac{2}{9}\m(G))^{1/3}$.
More precisely, he proved this for all hereditary families, not just families of all independent sets of a graph.
A family is {\it hereditary} if all subsets of its members are members.
Call a family $(r,t)$-\ekr\ if no $t$-intersecting subfamily of its $r$-sets is larger than some $t$-centered $r$-star.
Define $m(r,t)=(r-t)\binom{3r-2t-1}{t+1}+r$.

\begin{thm}
\label{t:Hered}
Let $\cF$ be any hereditary family of sets and let $t\le r$ be such that $m(r,t)\le\m(\cF)$.
Then, for any $\mt\not=S\ \sse\ \{t,\ldots,r\}$, the family $\cup_{s\in S}\cF^s$ is $(r,t)$-\ekr.
\end{thm}

The rough statement just mentioned is the instance $t=1$, $S=\{r\}$, and $(r-1)\binom{3(r-1)}{2}+r\le\m(\cF)$.
The theorem is not proven with the aid of shifting, but instead by the use of a version of the so-called Local LYM Inequality (see \cite{BollobasBook}, p.12) of Sperner \cite{Sperner} on sizes of shadows, generalized to hereditary families.

In \cite{BorgThesis}, Borg had already generalized Conjecture \ref{j:HolTal} to include all hereditary families.
In \cite{BorgExt}, he generalized it further, conjecturing the strengthening of Theorem \ref{t:Hered} that it holds with $m(r,t)=(t+1)(r-t+1)$, adding also that the conclusion holds strictly unless $\m(\cF)=m(r,t)$ and $S=\{r\}$.
He proved that this strengthening is true for shifted families $\cF$.

Also in 2009, Borg and Holroyd \cite{BorgHolrSingle} expanded upon Theorem \ref{t:Unions} significantly, as follows.
They call a graph a {\it modified} power of a cycle if it can be obtained from some $C_n^{(k)}$ by adding a set of consecutive edges of distance $k+1$ in $C_n$.

\begin{thm}
\label{t:MoreUnions}
Let $G$ be a disjoint union of complete graphs, copies of mnd graphs, powers of cycles, modified powers of cycles, and at least one singleton.
Then $G$ is $r$-\ekr\ for all $r\le \a(G)/2$.
\end{thm}

The list includes what they call {\it mnd graphs}, which are a generalization of path powers.
While vertices of a path are adjacent to all vertices to its right at distance at most some fixed $k$ to form a path power, the fixed $k$ condition is relaxed to a non-decreasing sequence of distances to form an mnd graph.
As graphs by themselves, Borg \cite{BorgEKRSep} later proved that mnd graphs satisfy Conjecture \ref{j:HolTal}.
Such graphs, of course, are chordal.
Observing that several verifications of Conjecture \ref{j:HolTal} involved chordal graphs and graphs with isolated vertices, the authors of \cite{HurlKamaChord} proved the following theorem in 2011.
Of course, as noted above, the role of an isolated vertex in these theorems is that such a vertex is a maximum center.

\begin{thm}
\label{t:Chordal}
If $G$ is a chordal graph with at least one singleton, then $G$ is $r$-\ekr\ for all $r\le \a(G)/2$.
\end{thm}

Disjoint unions of chordal graphs are chordal, and so it is not necessary to state the theorem with disjoint union language.
While stated only for chordal graphs, their proof actually shows that the statement is true for the larger class of dismantlable graphs.
By Dirac's theorem \cite{Dirac}, a graph is chordal if and only if it has a simplicial elimination ordering; that is, one can repeatedly remove a simplicial vertex --- one whose neighbors form a clique --- until only a complete graph remains.
Analogously, a graph is {\it dismantlable} if one can repeatedly remove a {\it corner} vertex --- one whose neighbors are neighbors of one of its neighbors --- until only a complete graph remains.%
\footnote{Nowakowski and Winkler \cite{NowaWink} showed that such graphs characterize cop-win graphs.}

More carefully, $v$ is a corner vertex if $N[v]\ \sse\ N[u]$ for some vertex $u\not=v$, where $N[x]=N(x)\cup\{x\}$ is the {\it closed} neighborhood of $x$; any such vertex $u$ is called a {\it witness}.
The proof of Theorem \ref{t:Chordal} uses an induction idea for vertices similar to Talbot's edge contraction operator used to prove Theorem \ref{t:Cycles}.
Here they define $G\da v=G-N[v]$ and note that both $G-v$ and $G\da v$ preserve dismantlability for all $v$.
Additionally, they observe that, if $v$ is a corner of $G$ and $u$ is a witness, then $\m(G-u)\ge\m(G)$ and $\m(G\da u)+1\ge\m(G)$.
These imply that, if $r\le\m(G)/2$, then $r\le\m(G-u)/2$ and $r-1\le\m(G\da u)/2$, which allows for induction to be applied.
From a dismantling order $(v_1,\ldots,v_n)$ on $V(G)$, they define a shift from any $v_i$ to $v_1$ and proceed with a simplified version of the partition ideas mentioned above, making use of the identity $|\cI_x^r(G)| = |\cI_x^r(G-v)|+|\cI_x^{r-1}(G\da v)|$ (for isolated $x$ and $v\not=x$) that drives the induction.

Soon after, Woodroofe \cite{WoodroofeSimp} proved a more general result for hereditary families via {\it (exterior) algebraic} shifting, as opposed to the {\it combinatorial} shifting discussed above.
Algebraic shifting was introduced by Kalai \cite{Kalai} and takes advantage of the geometry that arises from viewing a hereditary family as an abstract simplicial complex.
For example, the concepts of a simplicial elimination (or dismantling) ordering for chordal (or dismantlable) graphs generalize to shellability and the even more general sequentially Cohen-Macauley property for simplicial complexes, which we will not define here.
A graph $G$ is called {\it sequentially Cohen-Macaulay} if $\cI(G)$ is sequentially Cohen-Macaulay.
Among the graphs in this class are chordal graphs, graphs with no induced cycles of length other than 3 or 5, incomparability graphs of shellable posets, and others (see \cite{WoodroofeSimp}).
Woodroofe proved the following.

\begin{thm}
\label{t:CohMac}
If $G$ is a sequentially Cohen-Macauley graph with an isolated vertex, then $G$ is $r$-\ekr\ for all $r\le\m(G)$.
\end{thm}

Additionally, Woodroofe uses these tools to prove the most general version to date along the lines of Theorem \ref{t:MoreUnions}.

\begin{thm}
\label{t:MostUnions}
Let $G$ be the disjoint union of $k$ graphs, including a singleton.
Then $G$ is $r$-\ekr\ for all $r\le k/2$.
\end{thm}

Regarding graphs without singletons, in \cite{HurlKamaChord} we also find the verification of Conjecture \ref{j:HolTal} for a connected sequence of complete graphs called {\it chains}.
For a sequence of complete graphs $K_{a_1},\ldots,K_{a_m}$, we form the chain graph $G=G(a_1,\ldots,a_m)$ by allowing two consecutive $K_{a_i}$ and $K_{a_{i+1}}$ to share a single, distinct vertex, which requires of course that each $a_i\ge 2$.
The chain $G$ is called {\it special} if $m\le 2$ or if $m>2$, $a_i > a_{i-1}$ for each $2\le i<m$, and $a_m\ge a_{m-1}$.
They proved the following.

\begin{thm}
\label{t:Chains}
If $G$ is a disjoint union of at most two special chains, then $G$ is $r$-\ekr\ for all $r$.
\end{thm}

Disjoint unions $G$ of powers of cycles were considered by Hilton, Holroyd, and Spencer \cite{HilHolSpe}.
A {\it simple} cycle is just the cycle $C_n=C_n^{(1)}$.
For $r=\a(G)$, they proved that such unions are $r$-\ekr, provided that some component is simple, and were also able to identify which vertices are $r$-witnesses.
Included as ``cycles'' are the degenerate cases of a vertex ($C_1^{(1)}$) and an edge ($C_2^{(1)}$).
The problems for other values of $r$ remain open.
An $(a,b)$-{\it coloring} of a graph $G$ is a function 
$c:V(G)\rar\binom{[a]}{b}$ such that adjacent vertices receive disjoint subsets.
The {\it fractional chromatic number} of $G$, denoted $\xh$, is the minimum of $a/b$ over all $(a,b)$-colorings of $G$.

\begin{thm}
\label{t:CyclePowerUnions}
Let $G$ be a disjoint union of powers of cycles, at least one of which is simple.
Then $G$ is $\a(G)$-\ekr, and all $\a(G)$-witnesses belong to the component having least fractional chromatic number.
\end{thm}

Important ingredients in the proof of this result include the partitioning and shifting from Talbot's proof of Theorem \ref{t:Cycles}, as well as the fact that if $\xh(C_n^{(k)})\le \xh(C_{n'}^{(k')})$ then there is a homomorphism from $C_n^{(k)}$ to $C_{n'}^{(k')}$.
For the case of two simple cycles, they prove more.

\begin{thm}
\label{t:2Cycles}
If $G$ is a disjoint union of two cycles then $G$ is $r$-\ekr\ for all $r$; furthermore, all of the vertices
of any of its even cycles are $r$-witnesses.
\end{thm}

Hilton and Spencer \cite{HiltSpen,HiltSpenKing} then expanded on Theorem \ref{t:CyclePowerUnions}, handling all values of $r$.

\begin{thm}
\label{t:CyclePowerUnionsPath}
Let $G$ be a disjoint union of powers of cycles and a power of $H$, where $H$ is either a path or a cycle.
Write $G=H^k\cup\left(\cup_{i=1}^m C_{n_i}^{k_i}\right)$ and suppose that $k\le min_ik_i$.
Then $G$ is $\a(G)$-\ekr\ for all $r$, and some $r$-witness has minimum degree.
\end{thm}

Borg and Feghali \cite{BorgFegh} then prove the strictness of Theorem \ref{t:CyclePowerUnionsPath} when $k<\min_ik_i$.
Additionally, they give a simpler proof that relies on Katona's Shadow Intersection Theorem \cite{KatonaShadInt}.

In addition to studying the lexicographic product, as in Theorem \ref{t:LexProd}, it is natural to consider Cartesian products as well. 
For graphs $G$ and $H$, the {\it Cartesian product} $G\mathbin{\gbox} H$ has vertices $V(G)\times V(H)$ and edges $(u,x)(v,y)$ whenever $u=v$ and $xy\in E(H)$, or $uv\in E(G)$ and $x=y$.
For example, $K_2\mathbin{\gbox}K_2=C_4=Q^2$, and the $d$-{\it dimensional cube} $Q^d$ can be defined by $Q^d=K_2\mathbin{\gbox}Q^{d-1}$.
In \cite{HurlKamaLadders}, we find that the {\it ladder} $L_k=K_2\mathbin{\gbox}P_k$ is 3-\ekr.
Subsequently, Li and Zhang \cite{LiZhang} proved much more, using the traditional shifting and partitioning techniques outlined above.

\begin{thm}
\label{t:Ladders}
The ladder $L_k$ is $r$-\ekr\ for all $r$, and is strictly so unless $r=k-1$.
Additionally, the maximum centers are the minimum degree vertices of $L_k$.
\end{thm}


One of the simplest open cases for disjoint unions of graphs was the union of $k$ length-2 paths, denoted $kP_3$.
Generalizing Theorem \ref{t:Signed} (on $kK_2=kP_2$) in this manner, this case was partially resolved in 2020 \cite{FegHurKam} by also generalizing Bollob\'as and Leader's use of the idea behind Katona's cycle method.

\begin{thm}
\label{t:kP3}
The graph $G=kP_3$ is $r$-\ekr\ for all $r\le \m(G)/2=k/2$.
\end{thm}

Note that this also similar to Theorem \ref{t:Unions} without the isolated vertex.
The authors conjectured that $kP_3$ is $r$-\ekr\ for all $r$.
The main hurdle in generalizing Bollob\'as and Leader's proof is that, unlike $kP_2$, the graph $kP_3$ is not vertex transitive:
The $k$ degree-2 {\it middle} vertices live in one orbit and the $2k$ leaves live in another.
This causes the usual counting arguments to not line up properly, and the way around it is to state the cycle argument in a well-known probabilistic form, so that objects from the resulting probability space can be sampled in a controlled, but non-uniform manner that compensates for the imbalance.

A second hurdle arises in each component of $kP_3$ having the ability to contribute up to two vertices (leaf {\it siblings}) to an independent set, instead of up to one for $kP_2$.
For this, the authors still use ``good'' permutations of the leaves, but pair ``main'' intervals with opposite sub-intervals to distinguish leaves that appear alone or with their siblings in an independent set.
Further, the ``main'' intervals have an adjacent ``invisible'' interval that picks up middle vertices.
So that enumeration works well, this structure is reflected as well before rotating as in the original cycle method.

The obvious next generalization to attempt would be $kP_4$.
However, the authors suggest that unions of 3-claws also continue this sequence of graphs and might more easily succumb to a cycle-method-like symmetry argument --- a $k$-{\it claw} is the complete bipartite graph $K_{1,k}$, so $K_{1,1}=P_2$ and $K_{1,2}=P_3$.
Thus, an \ekr\ theorem for unions of claws would be similar to the relevant portion of Theorem \ref{t:Unions} without an isolated vertex.

An alternative to studying specifically structured families of graphs is to consider density; i.e. minimum or average degree.
For vertices $u$ and $v$ in a graph $G$, we use the notations $\deg_G(u)$, $\D(G)$, $\od(G)$, and $\dist_G(u,v)$ for the degree of vertex $u$, the maximum degree of $G$, the average degree of $G$, and the distance between $u$ and $v$ in $G$, respectively; we may omit the subscript if the context is clear.
In \cite{FranHurl} we find the following four theorems for sparse graphs.

\begin{thm}
\label{t:MaxDeg}
If $G$ is a graph on $n$ vertices with $\D(G)\le d$, then $G$ is $r$-\ekr\ for $r < (8n/27d)^{1/2}$. 
\end{thm}

Note that, for bounded-degree graphs, this improves the range of $r$ from the $O(n^{1/3})$ of Theorem \ref{t:Hered} to $O(n^{1/2})$.

\begin{thm}
\label{t:MaxAveDeg}
If $G$ is a graph on $n$ vertices with $\od(G)\le c$, for any $c\ge e/36$, then $G$ is $r$-\ekr\ for $r < (n/18c)^{1/3}$. 
\end{thm}

Here, bounded average graphs yield a range of the same order of magnitude as in the more general Theorem \ref{t:Hered}, but with an increased constant.

A {\it split} vertex is a vertex of degree at least 3.
A {\it spider} is a tree with exactly one split vertex.

\begin{thm}
\label{t:Spiders}
If $G$ is a spider on $n$ vertices, then $G$ is $r$-\ekr\ for $r\le (n\ln 2)^{1/2}-((\ln 2)/2)^{1/2}$.
\end{thm}

This returns the range of $r$ for spiders to $O(n^{1/2})$ and increases the constant found in Theorem \ref{t:MaxDeg}.

\begin{thm}
\label{t:Trees}
If $G$ is a tree on $n$ vertices, with exactly $s>1$ split vertices, then $G$ is $r$-\ekr\ for $2s < r \le (n\ln c)^{1/2}-((\ln c)/2)^{1/2}$, where $c=2-2s/r$.
\end{thm}

There is curiosity about there being a lower bound on $r$ here; the authors don't believe there should be one but didn't find a way to avoid it.

The key to these theorems is the method of {\it diversity} (originally called {\it unbalance} in \cite{LemoPalm}), which for a family $\cF$ is defined to be $\c(\cF) = |\cF|-|\cF_z|$ for any maximum center $z$.
The origins of the concept are found in a result of \cite{DinuFrie} that states all but at most $O(\binom{n-2}{k-2})$ members of any intersecting $\cF\ \sse \binom{[n]}{r}$ contain a common element.
The key lemma in \cite{FranHurl} is a result from \cite{FranklDiversity} that every intersecting $\cF\ \sse\binom{[n]}{r}$, with $n>72r$, has $\c(\cF)\le\binom{n-3}{r-2}$.
How this is used to prove Theorem \ref{t:MaxDeg}, for example, is to let $z$ be a maximum center of $\cF$ and choose some $E\in\cF-\cF_z$.
Because $G$ has bounded degree $d$, one can compute a simple lower bound $f(n,r,d)$ on the number of sets $F\in\cI_z^r(G)$ that are disjoint from $E$ and so can't be in $\cF$.
Then $|\cF| = |\cF_z|+\c(\cF) \le |\cI_z^r(G)| - f(n,r,d) + \binom{n-3}{r-2}$, which is at most $|\cI_z^r(G)|$ when $f(n,r,d)\ge \binom{n-3}{r-2}$.
Hence the bound on $r$ in the theorem comes from making the last inequality true.
The specific properties of the graphs in the other theorems are used to modify their bounds on $r$ similarly.

Another theorem on trees was given by Feghali, Hurlbert, and Kamat \cite{FegJohTho}.
A spider is $\ell$-{\it uniform} if every leaf has distance $\ell$ from the split vertex; it has also been called a {\it depth-$\ell$ claw} in the literature.
For example, the 1-uniform spider with $k$ leaves is $K_{1,k}$, which is $r$-\ekr\ for all $r\le k/2$ by Theorem \ref{t:EKR}.

\begin{thm}
\label{t:Claw2}
Every 2-uniform spider $G$ is $r$-\ekr\ for all $r\le (\m(G)+1)/2$.
In addition, if $G$ has $k$ leaves then $G$ is not $k$-\ekr.
\end{thm}

Finally, Holroyd and Talbot \cite{HolrTalb} opined that, if a counterexample $G$ to Conjecture \ref{j:HolTal} exists, then one exists with $\m(G)=\a(G)$.
This motivates studying such graphs, which are equivalent to the well-known class of well-covered graphs: 
a {\it vertex cover} is a set of vertices that touch all edges, and a graph is {\it well-covered} if all its minimal vertex covers have the same size.  
Examples of such graphs already studied above include disjoint unions of complete graphs, balanced complete bipartite graphs, and graphs of independence number 2 with no dominating vertex (e.g., complements of cycles).
A {\it pendant vertex} has degree 1.
Another construction of a well-covered graph is a {\it fully pendant graph} $G^*$, which is formed from a graph $G$ by adding to each of its vertices a pendant neighbor.
De Silva, Dionne, Dunkelberg, and Harris \cite{DeDiDuHa} recently proved the following.

\begin{thm}
\label{t:WellCovK}
The fully pendant graph $K_n^*$ is $r$-\ekr\ for all $r\le n/2$, and strictly so when $r<n/2$; all $r$-witnesses are pendant vertices.
\end{thm}


\section{\hk\ graphs}
\label{s:HK}

As in the case of isolated vertices, it can be part of a useful technique to know where the maximum $r$-centers are located in a graph, especially within inductive or injective methods.
In studying the \ekr\ property for connected graphs, it is reasonable to start with the fewest possible edges; i.e. trees.
Theorems \ref{t:Spiders}, \ref{t:Trees}, and \ref{t:Claw2} are important results in this direction.
As has been noted regarding the above results, among known witnesses are minimum degree vertices.
It makes sense, therefore, to guess that every tree has some leaf as a maximum center, as was conjectured in \cite{HurlKamaChord}, and verified for small $r$.

\begin{thm}
\label{t:SmallrHK}
If $T$ is a tree and $r\le 4$, then some maximum $r$-center of $T$ is a leaf.
\end{thm}

However, numerous counterexamples to the conjecture were soon found for $r\ge 5$ \cite{Baber,BorgStars,FegJohTho}.
The maximum star centers in \cite{Baber,BorgStars} were of degree 2.
However, Feghali, Johnson, and Thomas \cite{FegJohTho} showed that $r$-centers can be far from leaves and have any degree.

\begin{thm}
\label{t:DeepCenters}
For all $\ell\ge 3$ and $d\ge 2$ there is a tree $T$ and $r\le\a(T)$ with a maximum $r$-center of degree $d$ at distance $\ell$ from every leaf.
\end{thm}

Let $S^{t,\ell}$ be the $\ell$-uniform spider with split vertex of degree $t$.
Then define the tree $T=T^{d,t,\ell}$ by joining a vertex $v$ to the split vertices of $d$ copies of $S^{t,\ell-1}$.
Through very careful analysis, they are able to show that $|\cI_v(T)|>|\cI_u(T)|$ for every leaf $u$ when $t$ is large enough.
It follows, then, that $|\cI_v^r(T)|>|\cI_u^r(T)|$ for some value of $r$.
They complete the proof by showing that, for every remaining vertex $x$, $|\cI_v^r(T)|\ge |\cI_x^r(T)|$ for all $r$.

Based on these results, Estrugo and Pastine \cite{EstrPast} define a tree $T$ to be $r$-\hk\ if some leaf of $T$ is a maximum $r$-center, and simply \hk\ if it is $r$-\hk\ for all $r$.
More generally, one can define a graph $G$ to be $r$-\hk\ if some minimum-degree vertex of $G$ is a maximum $r$-center, and simply \hk\ if it is $r$-\hk\ for all $r$.
The problem then arises to discover which graphs are ($r$-)\hk.%
\footnote{Of course, regular graphs are!}

One of the first such results is found in \cite{HurlKamaTrees}, that spiders are \hk.
They not only show where the maximum $r$-center is, but also describe a partial order on the sizes of all $r$-centers.
Let $S=S(\ell_1,\ldots,\ell_k)$ denote the spider with $k$ leaves $v_1,\ldots,v_k$, with $v_i$ at distance $\ell_i$ from the split vertex $v$.
Say that $S$ is in {\it spider order} if the following three conditions hold for all $i<j$:
\begin{itemize}
    \item 
    if $\ell_i$ and $\ell_j$ are both odd then $\ell_i\le \ell_j$;
    \item 
    if $\ell_i$ and $\ell_j$ are both even then $\ell_i\ge \ell_j$; and
    \item 
    otherwise $\ell_i$ is odd and $\ell_j$ is even.
\end{itemize}

\begin{thm}
\label{t:SpiderHK}
Suppose that $S=S(\ell_1,\ldots,\ell_k)$ is in spider order, with split vertex $v$ and leaves $v_1,\ldots,v_k$ as defined above.
Then $S$ is \hk.
Moreover, for every $r$, $1\le i\le j\le k$, and degree-2 vertex $y$ on the $vv_i$-path in $S$ we have
\begin{enumerate}
    \item 
    $|\cI_y^r(S)|\le |\cI_{v_i}^r(S)|$,
    \item 
    $|\cI_v^r(S)|\le |\cI_{v_k}^r(S)|$, and
    \item 
    $|\cI_{v_j}^r(S)|\le |\cI_{v_i}^r(S)|$.
\end{enumerate}
In particular, $v_1$ is a maximum $r$-center.
\end{thm}

They proved each part of the theorem injectively, using a {\it flipping} operation on paths.
For example, if one flips the $yv_i$-path by mapping it onto itself in reverse order, independent $r$-sets containing $y$ get mapped to independent $r$-sets containing $v_i$, thereby injecting $\cI_y^r(S)$ into $\cI_{v_i}^r(S)$.
The same can be done for $v$ instead of $y$.
The interesting case is statement 3, of course, where the flipping operation is trickier, leading to an injection that works only for spider order.

In \cite{EstrPast}, Estrugo and Pastine generalize the flip to certain paths in graphs.
A {\it thread} in a graph $G$ is a path $P$ such that every interior vertex of $P$ has degree 2 in $G$.
Vacuously, an edge is a thread, for example.
They define a path $P=v_1\cdots v_h$ to be an {\it escape path} for $v_1$ if $v_1\cdots v_{h-1}$ is a thread and $v_h$ is pendant.
By flipping on $P$ they prove that $|\cI_{v_1}^r(G)|\le |\cI_{v_h}^r(G)|$, which yields the following theorem.

\begin{thm}
\label{t:Escape}
If every non-pendant vertex of a graph $G$ has an escape path, then $G$ is \hk.
\end{thm}

Hence caterpillars, for example, are \hk.
A graph $G$ is {\it pendant} if every split vertex of $G$ has a pendant neighbor.
Another example is that pendant graphs are \hk.
A pendant spider has also been called {\it short} spider.
For a vertex $x$ in $G$ one can {\it join a spider $S$ at} $x$ by identifying the split vertex of $S$ with $x$.
Then, for a set $D$ of vertices of $G$, a {\it short $D$-spidering of} $G$ joins to each vertex of $D$ a short spider.
It follows that if $D$ is a dominating set of $G$ then any short $D$-spidering of $G$ is \hk.
Define a set of vertices $D$ to be a {\it thread dominating set} of $G$ if every vertex not in $D$ has a thread to some vertex in $D$.
Estrugo and Pastine then observe the following, more general consequence of Theorem \ref{t:Escape}.

\begin{thm}
\label{t:ShortSpidering}
If $D$ is a thread dominating set of $G$, then any short $D$-spidering of $G$ is \hk.
\end{thm}

Hence all pendant graphs are \hk.

The authors of \cite{HurlKamaTrees} conjectured that trees with no degree-2 vertices are \hk.
As a test case, the authors of \cite{FranHurl} propose showing that the following family of graphs is $r$-\ekr.
For $1\le i\le 3$, let $T_i(h)$ be a complete binary tree of depth $h$, having unique degree-2 vertex $v_i$.
Now join a new vertex $v$ to each $v_i$ to create the tree $T(h)$, in which every interior vertex has degree 3.

Finally, Zhang, Ran, and Huang \cite{ZhaRanHua} prove the following result, similar to Theorem \ref{t:SmallrHK}.
A {\it unicyclic graph} contains exactly one cycle.

\begin{thm}
\label{t:Uni}
If $G$ is a connected, unicyclic graph with at least one pendant vertex then $G$ is $r$-\hk\ for $r\le 4$.
There exist such graphs that are not $r$-\hk\ for $r\ge 5$.
\end{thm}

They also construct a caterpillar $P^{k,t}$ by attaching $t$ pendant vertices to each vertex of $P_k=v_1\cdots v_k$; thus $P^{k,1}=P_k^*$, the fully pendant contruction for well-covered graphs.
They prove, for all $r$ and $t$, that every pendant neighbor of $v_2$ or $v_{k-1}$ is a maximum $r$-center (the $t=1$ case had been done in \cite{DeDiDuHa}).

Zhang, et al. finish by offering the following conjecture.

\begin{cnj}
\label{j:LeafCenterNbr}
If $T$ is a caterpillar, $r\ge 1$, and $z$ is a maximum $r$-center, then the neighbor of $z$ has maximum degree.
\end{cnj}


\section{Chv\'atal's Conjecture}

One cannot avoid Chv\'atal's Conjecture \cite{Chvatal} when discussing this topic.
It is the elephant in the room.
We say that a family $\cF$ of sets is \ekr\ if some element $x$ satisfies $|\cF_x|\ge |\cH|$ for every intersecting family $\cH\ \sse\ \cF$.

\begin{cnj}
\label{j:Chvatal}
{\bf (Chv\'atal, 1972)}
Every hereditary family $\cF$ is \ekr.
\end{cnj}

It is a monument to the difficulty of the conjecture that it is so well-known and yet so few papers have addressed it in the half century since its birth.
There were some early papers.
Kleitman and Magnanti \cite{KleiMagn} solved the case in which $\cF$ is the union of two stars.
Sterboul \cite{Sterboul} solved the rank at most 3 case, where the {\it rank} of $\cF$ is the size of its largest set.
In 1976, Sch\"onheim \cite{Schonheim} solved case in which the maximal elements of $\cF$ intersect.
Then there was a long gap until Snevily \cite{Snevily} solved the case in which $\cF$ is 1-shifted in 1992, and quite another long period before Borg \cite{BorgChvatal} generalized the conjecture to weighted families in 2011, solving the 1-shifted case as well.
Then in 2017, the paper \cite{CzaHurKam} reproved Sterboul's result, once with a relatively similar approach, and once with the Sunflower Lemma, in hopes that the method might yield more fruit.
Subsequently, Olarte, Santos, and Spreer \cite{OlaSanSpr} reproved both Kleitman-Magnanti and Sterboul with shorter proofs.
Most recently, Frankl and Kupavskii \cite{FranKupaMatch} solved the case in which $\cF$ has cover number 2.

{
\color{brwn}
\noindent
\hrulefill
\medskip

\centerline{
***
$\quad$
This is not stated correctly.  The correct theorem of Frankl and Kupavskii is the following.
$\quad$
***
}
\medskip

{\it 
\noindent
Let $\cF$ be hereditary and $\cH\sse\cF$ be intersecting.
If $\t(\cH)\le 2$ then there is some $x$ such that $|\cF_x|\ge |\cH|$.
}
\medskip

\noindent
\hrulefill
}
\medskip

In the meantime, some of the results above actually validate Conjecture \ref{j:Chvatal} when $\cF=\cI(G)$ for various graphs $G$ that are $r$-\ekr\ for all $r$.
For example, in Theorem \ref{t:Cycles}, because a cycle power $G=C_n^{(k)}$ is vertex transitive, every vertex is a maximum $r$-center for every $r$.
This means that, for any fixed vertex $x$, if $\cH\ \sse\ \cI(G)$, then $|\cH| = \sum_r|\cH^r| \le \sum_r|\cI_x^r(G)| = |\cI_x(G)|$.
This not only works for vertex transitive graphs, but for graphs for which some vertex is an $r$-witness for all $r$: e.g., path powers $P_n^{(k)}$ (Theorem \ref{t:Paths}, minimum degree vertex), lexicographic products $P_n^{(k)}[K_m]$ (Theorem \ref{t:LexProd}, minimum degree vertex), disjoint union of two cycles, one of which is even (Theorem \ref{t:2Cycles}, even cycle vertex), and ladders (Theorem \ref{t:Ladders}, minimum degree vertex).
This highlights the importance of tracking maximum $r$-centers.


\section{Questions and Remarks}

This subject is clearly rich and vibrant, teeming with tantalizing questions, problems, variations, and conjectures.
Surely the reader will think of lovely problems to attack beyond what is included below, but here are a few thoughts to try in hopes of making progress.

\begin{enumerate}
\item 
\ekr\ problems.
    \begin{enumerate}
    \item 
    What Cartesian (or other) products of various graphs are $r$-\ekr; e.g., products of paths, cycles, or complete graphs?
    For some graph product $\times$, could it be that, if $G$ and $H$ are $r$-\ekr\ for all $r$, then is $G\times H$ $r$-\ekr\ for all $r$?
    \item 
    Are other well-covered graphs $r$-\ekr?
    The {\it rook graph} $K_m\mathbin{\gbox} K_m$ is particularly attractive because it is both well-covered and a Cartesian product.
    \item 
    Might the techniques applied to sparse graphs be tweaked somewhat so as to work for dense graphs?
    \item 
    Are all spiders $r$-\ekr?
    Is the $\ell$-uniform spider with $k$ leaves $r$-\ekr\ if and only if $kP_{\ell+1}$ is $r$-\ekr?
    \end{enumerate}
\item 
\hk\ problems.
    \begin{enumerate}
    \item 
    Is the internally 3-regular tree $T(h)$ $r$-\hk\ for all $r$?
    \item 
    Where are the maximum $r$-centers for caterpillars?
    (See \cite{EstrPast}.)
    \item 
    What other graphs $G$ (besides cycles) with $\d(G)=2$ are $r$-\hk?
    \end{enumerate}
\end{enumerate}

\section*{Acknowledgement}

This paper arose from a talk I presented at the July, 2024, Summit280 Conference in Budapest in celebration of the $70^{\rm th}$ birthdays of P\'eter Frankl, Zolt\'an F\"uredi, Ervin Gy\H ori, and J\'anos Pach, two of whom were my professors in graduate school and had a strong influence in my pursuit of combinatorial mathematics.
The author extends his gratitude to the organizers of the conference for inviting me to write this survey.


\bibliographystyle{acm}
\bibliography{holtal}

\end{document}